\documentclass{segabs}

\usepackage{amsmath}
\usepackage{amsfonts}
\usepackage{bm} 

\usepackage[pass,paperwidth=8.5in,paperheight=11in]{geometry}

\newtheorem{alg}{Algorithm}

\begin{document}

\title{Nonlinear seismic imaging via reduced order model backprojection}

\renewcommand{\thefootnote}{\fnsymbol{footnote}} 

\author{Alexander V. Mamonov\footnotemark[1], University of Houston;
Vladimir Druskin and Mikhail Zaslavsky, Schlumberger}

\righthead{Nonlinear imaging via ROM backprojection}

\maketitle

\section{Summary}

We introduce a novel nonlinear seismic imaging method based on 
model order reduction. The reduced order model (ROM) is an orthogonal
projection of the wave equation propagator operator on the subspace
of the snapshots of the solutions of the wave equation. It can be 
computed entirely from the knowledge of the measured time domain 
seismic data. The image is a backprojection of the ROM using the 
subspace basis for the known smooth kinematic velocity model. The 
implicit orthogonalization of solution snapshots is a nonlinear 
procedure that differentiates our approach from the conventional 
linear methods (Kirchhoff, RTM). It allows for the removal of 
multiple reflection artifacts. It also enables us to estimate the 
magnitude of the reflectors similarly to the true amplitude 
migration algorithms.

\section{Introduction}

To simplify the exposition we make a number of assumptions on our model.

First, we consider an acoustic wave equation
\begin{equation}
u_{tt} = \tensor{A} u, \quad t \in [0, T],
\label{eqn:waveeqn}
\end{equation}
for the pressure $u$. 
We treat the spatial operator 
$\tensor{A} = \tensor{C}^2 \bm\Delta$ as a matrix 
$\tensor{A} \in \mathbb{R}^{N \times N}$, a discretization on some fine grid.
The diagonal matrix $\tensor{C} = \mbox{diag}(\mathbf{c})$ contains the acoustic 
velocity values $c \in \mathbb{R}^N$ at the $N$ nodes of the fine grid, 
while $\bm\Delta$ discretizes the Laplacian on that grid.

Second, we assume that the sources can be modeled as an initial condition
\begin{equation}
\mathbf{u}(0) = \tensor{S}, \quad \mathbf{u}_t(0) = \mathbf{0},
\label{eqn:initcond}
\end{equation}
which can be easily achieved by considering the even part of the time 
domain solution of equation \ref{eqn:waveeqn} and thus the even part of the data.
The matrix $\tensor{S} \in \mathbb{R}^{N \times p}$ contains all $p$ sources
and is localized near the surface. 

Under the initial conditions of equation \ref{eqn:initcond} the solution to 
equation \ref{eqn:waveeqn} takes the form
\begin{equation}
\mathbf{u}(t) = \cos(t \sqrt{-\tensor{A}}) \tensor{S},
\label{eqn:solu}
\end{equation}
where the cosine and the square root are understood as matrix functions.
Note that the matrix function of time 
$\mathbf{u}(t) \in \mathbb{R}^{N \times p}$ consists of $p$ columns that 
contain the solutions corresponding to each source in $\tensor{S}$.

Third, we suppose that the source matrix $\tensor{S}$ admits a 
representation
\begin{equation}
\tensor{S} = q^2(\tensor{A}) \tensor{C} \tensor{E},
\label{eqn:src}
\end{equation}
where $\tensor{E} \in \mathbb{R}^{N \times p}$ are $p$ point sources 
supported on the surface and $q^2(\omega)$ is the Fourier transform of 
the source wavelet. Here we consider a Gaussian source wavelet
\begin{equation}
q^2(\tensor{A}) = e^{\sigma \tensor{A}}.
\end{equation}
For small $\sigma$ all quantities in equation \ref{eqn:src} are 
localized near the surface. We assume that the velocity near the 
surface is known, thus we know $\tensor{S}$.

Fourth, we assume that the sources and receivers are collocated. This 
assumption makes the construction of the reduced order model in the next 
section a lot more straightforward. However, our approach can be generalized 
to the much more realistic case of noncollocated sources and receivers. 
The particular form of the receiver matrix 
$\tensor{R} \in \mathbb{R}^{N \times p}$ that we use is
\begin{equation}
\tensor{R} = \tensor{C}^{-1} \tensor{E}.
\label{eqn:recv}
\end{equation}

Finally, under all the assumptions above we can write a data model
\begin{equation}
\tensor{F}(t) = \tensor{R}^T \cos(t \sqrt{- \tensor{A}}) \tensor{S},
\label{eqn:datamodel}
\end{equation}
which is a $p \times p$ matrix function of time.

The problem of seismic imaging that we solve here is to find an estimate
of the unknown velocity $\mathbf{c}$ from the knowledge of the time domain
data $\tensor{F}(t)$ for $t \in [0, T]$ and a smooth kinematic velocity
model $\mathbf{c}_0$.

\section{Time-domain interpolatory reduced order model}

At the core of our approach is the construction of the ROM 
that interpolates the measured seismic data.
The use of model order reduction techniques in
inversion was proposed in \cite{druskin2013solution, borcea2014model}
for parabolic (controlled source electromagnetic method, CSEM) inverse
problems. Unlike the diffusive parabolic case, where the authors 
employed frequency (Laplace) domain interpolation, the appropriate setting
for the wave equation inversion is the time domain. 

The particular form of the source and receiver matrices in equations
\ref{eqn:src}--\ref{eqn:recv} allows us to rewrite the data model 
from equation \ref{eqn:datamodel} in the completely symmetric form
\begin{equation}
\tensor{F}(t) = \widehat{\tensor{B}}^T 
\cos \left( t \sqrt{-\widehat{\tensor{A}}} \right) \widehat{\tensor{B}},
\label{eqn:symdata}
\end{equation}
where the symmetrized spatial operator is 
$\widehat{\tensor{A}} = \tensor{C} \bm\Delta \tensor{C}$ and the 
source/receiver matrix is given by 
$\widehat{\tensor{B}} = q(\widehat{\tensor{A}}) \tensor{E}$. 
Here we used the fact that analytic matrix functions commute with 
similarity transforms
and also that the symmetric operator $\widehat{\tensor{A}}$ is similar 
to the original $\tensor{A}$ with a similarity transform $\tensor{C}$.

In practice, the time domain data is measured at discrete time instants 
that we denote by $t_k = k \tau$ with $k = 0,1,\ldots,2n-1$, where 
$\tau$ is the sampling interval and $t_{2n-1} = T$ is the terminal time.

The discrete data samples $\tensor{F}_k = \tensor{F}(t_k)$ admit a 
representation
\begin{equation}
\tensor{F}_k = \tensor{F}(k \tau) = 
\widehat{\tensor{B}}^T \cos \left( k \; \mbox{arccos} 
\cos \tau \sqrt{-\widehat{\tensor{A}}} \right) \widehat{\tensor{B}} =
\widehat{\tensor{B}}^T T_k(\widehat{\tensor{P}}) \widehat{\tensor{B}},
\end{equation}
where $T_k$ are Chebyshev polynomials of the first kind and 
$\widehat{\tensor{P}}$ is the propagator given by 
\begin{equation}
\widehat{\tensor{P}} = \cos \left( \tau \sqrt{-\widehat{\tensor{A}}} \right).
\label{eqn:prop}
\end{equation}

We wish to construct a reduced order model of size $np \ll N$ that 
matches the $2n$ measured data samples exactly
\begin{equation}
\tensor{F}_k = 
\widehat{\tensor{B}}^T T_k(\widehat{\tensor{P}}) \widehat{\tensor{B}} = 
\widetilde{\tensor{B}}^T T_k(\widetilde{\tensor{P}}) \widetilde{\tensor{B}},
\quad k=0,1,\ldots,2n-1,
\label{eqn:interp}
\end{equation}
where $\widetilde{\tensor{P}} \in \mathbb{R}^{np \times np}$ and 
$\widetilde{\tensor{B}} \in \mathbb{R}^{np \times p}$. Since we are 
solving the inverse problem of seismic imaging the ROM 
$(\widetilde{\tensor{P}}, \widetilde{\tensor{B}})$ should be 
computable from the knowledge of the sampled data $\tensor{F}_k$ only.

It appears that the solution to the data interpolation problem 
of equation \ref{eqn:interp} can be found in the projection form
\begin{equation}
\widetilde{\tensor{P}} = \tensor{V}^T \widehat{\tensor{P}} \tensor{V},
\quad \widetilde{\tensor{B}} = \tensor{V}^T \widehat{\tensor{B}},
\label{eqn:proj}
\end{equation}
where the columns of $\tensor{V} \in \mathbb{B}^{N \times np}$ 
constitute an orthonormal basis for the subspace spanned by the 
discrete time snapshots of solutions
\begin{equation}
\widehat{\mathbf{u}}_k = \widehat{\mathbf{u}}(t_k) = 
\cos \left( k \tau \sqrt{- \widehat{\tensor{A}}} \right) \widehat{\tensor{B}} =
T_k(\widehat{\tensor{P}}) \widehat{\tensor{B}} \in \mathbb{R}^{N \times p}.
\label{eqn:snapshot}
\end{equation}
If we introduce the matrix of solution snapshots
\begin{equation}
\tensor{U} = \left[ \widehat{\mathbf{u}}_0, \widehat{\mathbf{u}}_1, \ldots,
\widehat{\mathbf{u}}_{n-1} \right] \in \mathbb{R}^{N \times np},
\end{equation}
then $\tensor{V}$ is defined simply by
\begin{equation}
\mbox{colspan} \;\tensor{V} = \mbox{colspan} \;\tensor{U}, \quad 
\tensor{V}^T \tensor{V} = \tensor{I}.
\end{equation}
Note that the above definition is not unique as $\tensor{V}$ is defined
up to an orthonormal change of variables in the projection subspace.
For the purpose of imaging some choices of $\tensor{V}$ are better
than others.

The orthogonalization of snapshots $\tensor{U}$ must respect the 
causality and the propagating nature of the time domain solutions of 
the wave equation. 
Thus, each snapshot should be orthogonalized only against the previous ones. 
In linear algebra this is known as Gram-Schmidt procedure or the QR
decomposition. Since our snapshots are matrices with $p$ columns
corresponding to all sources/receivers, we need a block version of
QR decomposition
\begin{equation}
\tensor{U} = \tensor{V} \tensor{L}^T,
\label{eqn:blockqr}
\end{equation}
where $\tensor{L}^T \in \mathbb{R}^{np \times np}$ is block upper 
triangular with blocks of size $p$.

Obviously, we cannot simply use equation \ref{eqn:blockqr}
since the snapshots $\tensor{U}$ are unknown to us. However, from 
the data we can obtain the inner products between the snapshots.
A basic multiplication property of Chebyshev polynomials
\begin{equation}
T_i(x) T_j(x) = \frac{1}{2} ( T_{i+j}(x) + T_{|i-j|}(x) )
\label{eqn:chebmult}
\end{equation}
combined with equation \ref{eqn:snapshot} immediately implies that
\begin{equation}
(\tensor{U}^T \tensor{U})_{ij} = 
\widehat{\mathbf{u}}_i^T \widehat{\mathbf{u}}_j =
\frac{1}{2}( \tensor{F}_{i+j} + \tensor{F}_{|i-j|} ).
\label{eqn:utu}
\end{equation}
Applying equation \ref{eqn:chebmult} twice we can also obtain
\begin{equation}
\begin{split}
&(\tensor{U}^T \widehat{\tensor{P}} \tensor{U})_{ij} = 
\widehat{\mathbf{u}}_i^T \widehat{\tensor{P}} \widehat{\mathbf{u}}_j = \\
& = \frac{1}{4}( \tensor{F}_{i+j+1} + \tensor{F}_{|i-j+1|} + 
\tensor{F}_{|i+j-1|} + \tensor{F}_{|i-j-1|} ).
\end{split}
\label{eqn:utpu}
\end{equation}

The knowledge of Gram matrix $\tensor{U}^T \tensor{U}$ from 
equation \ref{eqn:utu} allows us to compute the block lower 
triangular factor $\tensor{L}$ in equation \ref{eqn:blockqr}
via a block Cholesky decomposition
\begin{equation}
\tensor{U}^T \tensor{U} = \tensor{L} \tensor{L}^T.
\label{eqn:blockchol}
\end{equation}

Once the Cholesky factor is known we use equation \ref{eqn:utpu} to
obtain the final expression for the ROM
\begin{equation}
\widetilde{\tensor{P}} = 
\tensor{L}^{-1} (\tensor{U}^T \widehat{\tensor{P}} \tensor{U}) \tensor{L}^{-T},
\label{eqn:prom}
\end{equation}
entirely from the sampled data $\tensor{F}_k$.

\section{Backprojection imaging}

After the reduced order model of equation \ref{eqn:prom} is obtained 
from the measured data we need to extract from it the information about 
the velocity $\mathbf{c}$. The first step is to go from the ROM for the 
propagator $\widehat{\tensor{P}}$ to the reduced model for the 
symmetrized spatial operator 
$\widehat{\tensor{A}} = \tensor{C} \bm\Delta \tensor{C}$ by 
approximately inverting equation \ref{eqn:prop} using the first two 
terms in the Taylor's expansion
\begin{equation}
\widetilde{\tensor{A}} = \frac{2}{\tau^2} ( \widetilde{\tensor{P}} - \tensor{I} )
\approx \tensor{V}^T \widehat{\tensor{A}} \tensor{V}.
\label{eqn:arom}
\end{equation}

There are multiple ways to obtain the estimate of the velocity from
the knowledge of $\widetilde{\tensor{A}}$. One may employ optimization
to solve for $\mathbf{c}$ by minimizing the ROM misfit. Such procedure
is superior to the conventional full waveform inversion (FWI) which
minimizes the data misfit. However, in this work we are interested in
a non-iterative imaging algorithm that assumes the knowledge of a 
smooth kinematic background model denoted by $\mathbf{c}_0$.

If the projection subspace $\mbox{colspan} \; \tensor{V}$ is sufficiently
rich, then the backprojection must be a good approximation of the 
spatial operator
\begin{equation}
\widehat{\tensor{A}} \approx \tensor{V} \widetilde{\tensor{A}} \tensor{V}^T.
\label{eqn:bpapprox}
\end{equation}
However, we have no direct access to the orthonormal basis $\tensor{V}$.
We approximate it with a known basis $\tensor{V}_0$ for the smooth kinematic 
velocity model $\mathbf{c}_0$.

In order to get an imaging formula we also notice that the diagonal
of $\widehat{\tensor{A}}$ is proportional to the square of the velocity
\begin{equation}
\mathbf{c}^2 \propto \mbox{diag}(\widehat{\tensor{A}}) = 
\mbox{diag}(\tensor{C} \bm\Delta \tensor{C}),
\end{equation}
where the square $\mathbf{c}^2$ is understood componentwise. Similarly,
for the difference between the unknown velocity and the kinematic model
we can write
\begin{equation}
\delta \mathbf{c}^2 = \mathbf{c}^2 - \mathbf{c}^2_0
\propto \mbox{diag}(\widehat{\tensor{A}} - \widehat{\tensor{A}}_0).
\label{eqn:deltac2exact}
\end{equation}

Replacing the symmetrized operators $\widehat{\tensor{A}}$ and
$\widehat{\tensor{A}}_0$ in equation \ref{eqn:deltac2exact} with 
their backprojection approximations from equation \ref{eqn:bpapprox} 
and also using the approximation $\tensor{V} \approx \tensor{V}_0$ we
arrive at the formula
\begin{equation}
\delta \mathbf{c}^2 \propto \mbox{diag} 
\left( \tensor{V}_0 ( \widetilde{\tensor{A}} - \widetilde{\tensor{A}}_0 ) \tensor{V}_0^T \right),
\label{eqn:formula}
\end{equation}
where $\widetilde{\tensor{A}}$ is computed from the measured data
and $\widetilde{\tensor{A}}_0$, $\tensor{V}_0$ are easily found 
since $\mathbf{c}^2_0$ is known.

Observe that unlike the conventional imaging approaches (reverse time migration, 
Kirchhoff migration) the formula in equation \ref{eqn:formula} is nonlinear
in the measured data. This is due to the nonlinearity of the block Cholesky 
decomposition in equation \ref{eqn:blockchol} and inversion of the 
Cholesky factor $\tensor{L}$ in equation \ref{eqn:prom}. The nonlinearity 
that amounts to the implicit orthogonalization of solution snapshots 
$\tensor{U}$ allows for a better quality image. In particular, the 
orthogonalization removes the multiple reflection artifacts which 
are otherwise very difficult to deal with using conventional linear 
migration algorithms. This is illustrated in 
Figure~\ref{fig:ctrue_bproj1ordb_layers_n12,img_bproj1ordb_layers_n12,grad_bproj1ordb_layers_n12}
where we show a simple synthetic model with two layers. 
The backprojection images the layers correctly while suppressing 
the multiple reflection artifacts that are present in the RTM image as
ghost layers below the actual ones.

\multiplot{3}{ctrue_bproj1ordb_layers_n12,img_bproj1ordb_layers_n12,grad_bproj1ordb_layers_n12}{width=0.37\textwidth}
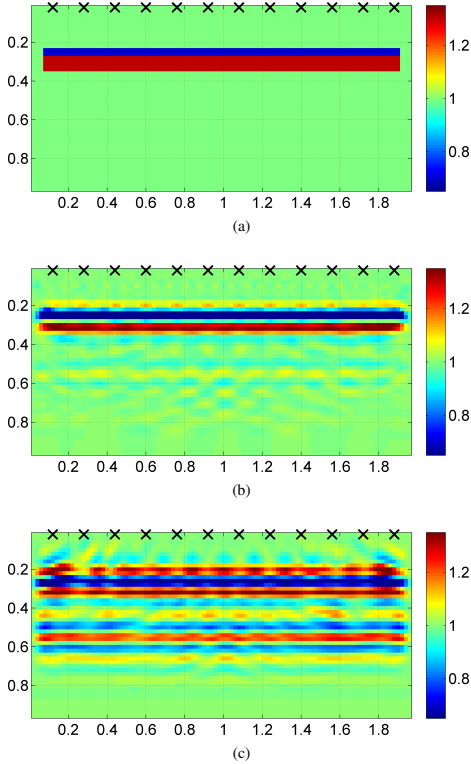
{Removal of multiples: (a) true velocity $\mathbf{c}$; 
(b) backprojection image $\mathbf{c}^\star$;
(c) RTM image computed as a gradient of conventional FWI.
Distances are in $km$, velocities in $km/s$, $p=12$ sources/receivers
are black $\times$.}

Equation \ref{eqn:formula} can be used in a number of ways to obtain the
image. The ambiguity comes from the choice of the proportionality factor. 
Here we choose it to be the background velocity, which leads us to a 
multiplicative imaging formula
\begin{equation}
\mathbf{c}^\star = \mathbf{c}_0 \sqrt{1 + \alpha \delta \mathbf{c}^2},
\label{eqn:imaging}
\end{equation}
where $\mathbf{c}^\star \approx \mathbf{c}$ is the image, $\alpha$ is a scalar
step length and all algebraic operations on the right hand side are performed
componentwise. With the imaging formula from equation \ref{eqn:imaging} 
at hand we can summarize our seismic imaging method in the following algorithm.

\begin{alg}[Nonlinear ROM backprojection imaging]~\\
\vskip-0.6in
~\\
\begin{enumerate}
\item Choose the sampling time interval $\tau$ and measure 
the discrete time samples of the seismic data with
equation \ref{eqn:datamodel}: 
$\tensor{F}_k = \tensor{F}(\tau k)$ for $k=1,2,\ldots,2n-1$.
\item Using equation \ref{eqn:utu} compute the snapshot 
Gram matrix $\tensor{U}^T \tensor{U}$ from the data and 
perform its block Cholesky decomposition 
as in equation \ref{eqn:blockchol}.
\item Using equation \ref{eqn:utpu} compute the matrix 
$\tensor{U}^T \widehat{\tensor{P}} \tensor{U}$ from the data 
and use Cholesky factor from step 2 to form the ROM 
$\widetilde{\tensor{P}}$ with equation \ref{eqn:prom}.
\item Choose a smooth kinematic velocity model $\mathbf{c}_0$
and use the block QR decomposition of equation \ref{eqn:blockqr}
to compute the orthonormal basis $\tensor{V}_0$;
project $\widehat{\tensor{P}}_0$ on $\tensor{V}_0$ 
using equation \ref{eqn:proj} to obtain the kinematic 
model ROM $\widetilde{\tensor{P}}_0$.
\item From the propagator ROMs $\widetilde{\tensor{P}}$ and 
$\widetilde{\tensor{P}}_0$ obtain the spatial operator ROMs 
$\widetilde{\tensor{A}}_0$ and $\widetilde{\tensor{A}}_0$ 
using equation \ref{eqn:arom}.
\item With the operator ROMs from step 5 and the orthonormal 
basis from step 4 compute $\delta \mathbf{c}^2$ from 
equation \ref{eqn:formula} and use the imaging formula
of equation \ref{eqn:imaging} to form the final image 
$\mathbf{c}^\star$.
\end{enumerate}
\end{alg}

\section{Numerical example: Marmousi model}

We evaluate the performance of our method on the synthetic data computed
for the Marmousi model by \cite{bourgeois1991marmousi}. 

\plot{log_full_ind460_bproj1ordb_marmfull_n10_m34-crop}{width=0.45\textwidth}
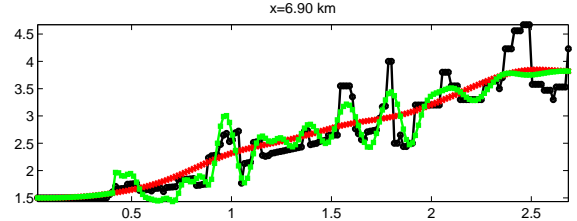
{Vertical well log for the Marmousi model at offset $x=6.9 \; km$.
Depth (horizontal axis) is in $km$, velocity (vertical axis) is in $km/s$. 
True $\mathbf{c}$ is black $\circ$, 
smooth kinematic velocity model $\mathbf{c}_0$ is red $\times$
and the image $\mathbf{c}^\star$ is green $\square$.}

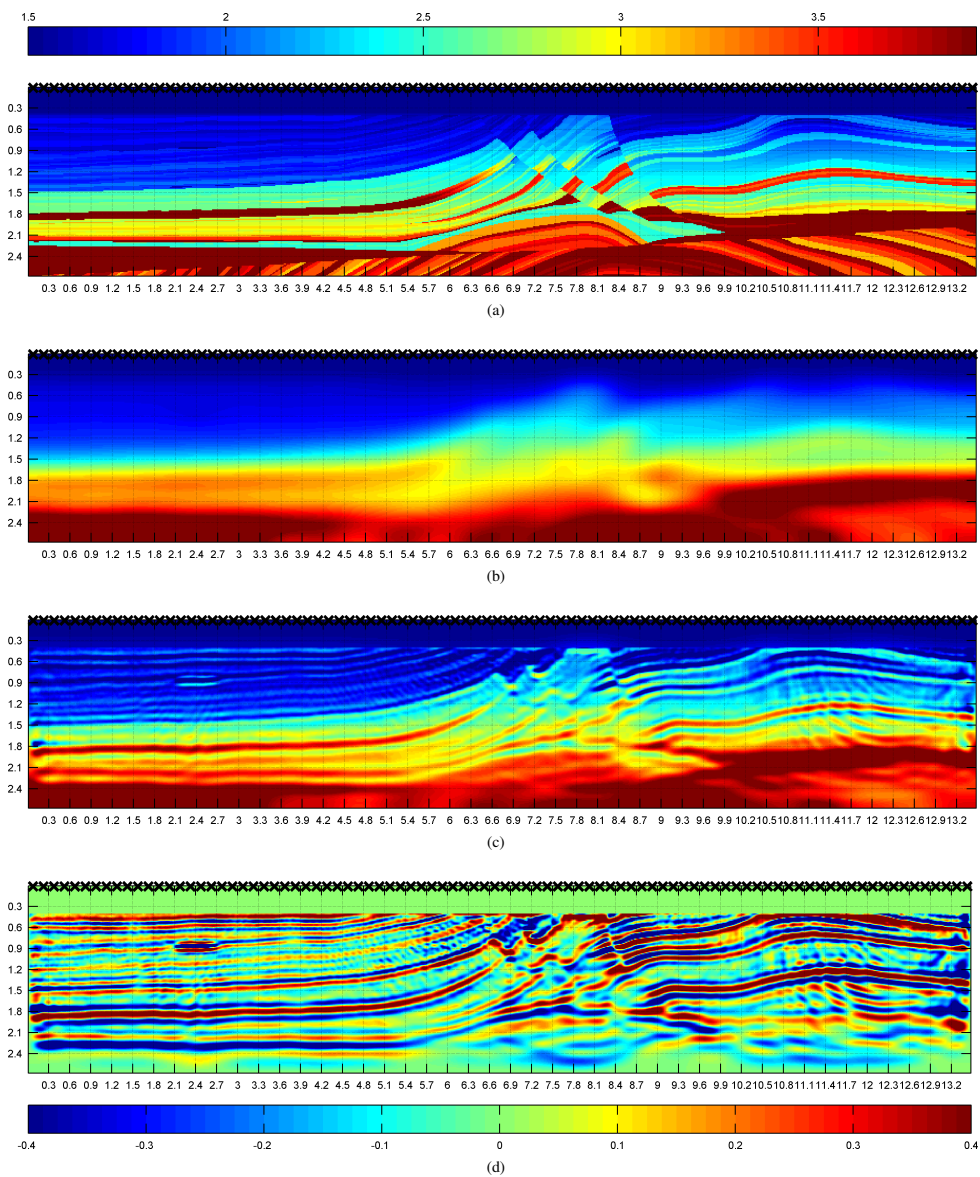
\multiplot*{4}{ctrue_bproj1ordb_marmfull_n10,c0_bproj1ordb_marmfull_n10,img_bproj1ordb_marmfull_n10,ddbp_bproj1ordb_marmfull_n10}{width=0.78\textwidth}{Seismic image for a Marmousi model: 
(a) True velocity $\mathbf{c}$; (b) Smooth kinematic velocity model $\mathbf{c}_0$; 
(c) Nonlinear ROM backprojection image $\mathbf{c}^\star$; 
(d) Difference $\mathbf{c}^\star - \mathbf{c}_0$. 
All distances are in $km$, velocities in $km/s$. 
The sources/receivers are black $\times$.}

The model is on a $15 \; m$ grid with $N = 900 \times 180 =$ $162,000$ nodes. 
The choice of the time interval $\tau$ is very important for our method. 
To make the orthogonalization procedure well conditioned, it should be chosen
at a Nyquist rate for the given source wavelet. Here we use $\tau = 33.5 \; ms$
($n=35$) which corresponds to the frequency of about $\omega = 15 \; Hz$. 
The data is measured for $p=90$ sources/receivers spaced uniformly every 
$150 \; m$. The kinematic model $\mathbf{c}_0$ is obtained by convolving 
the true velocity $\mathbf{c}$ with a Gaussian kernel of width $465 \; m$
and height of $315 \; m$.

In Figure~\ref{fig:ctrue_bproj1ordb_marmfull_n10,c0_bproj1ordb_marmfull_n10,img_bproj1ordb_marmfull_n10,ddbp_bproj1ordb_marmfull_n10} we show the image $\mathbf{c}^\star$ and the difference
$\mathbf{c}^\star - \mathbf{c}_0$ along with the smooth kinematic and true
Marmousi models. We observe very good recovery of all the model's features
down to $2.4 \; km$. The very bottom is not imaged because we had to truncate 
the data sampling at $2n = 70$ to avoid the reflections from the bottom, which
at the moment employs reflective boundary conditions instead of the PML.

We also show in Figure~\ref{fig:log_full_ind460_bproj1ordb_marmfull_n10_m34-crop}
a vertical well log. It demonstrates that our method performs well not only 
recovering the locations of the reflectors but also their strengths. 
We observe that the magnitude of the imaged velocity $\mathbf{c}^\star$ is in 
good agreement with the true model $\mathbf{c}$. In this particular aspect
our algorithm performs similarly to the true amplitude migration methods.

\section{Conclusions and discussion}

We introduced a novel nonlinear seismic imaging method based on the backprojection 
of reduced order models computed directly from the measured time domain data. 
The results of the early numerical experiments with synthetic data for 
Marmousi model show great promise.

The main issue to be solved to make the method viable for the real field data
is to remove the assumption that the sources and receivers are collocated.
This is certainly possible if one uses different left (source) and right 
(receiver) subspaces in the projection equation \ref{eqn:proj}. 
Other possible improvements include the implementation of absorbing boundary
conditions (PML) and more accurate source models.




\bibliographystyle{seg}  
\bibliography{biblio}

\end{document}